\documentclass{amsart}
\usepackage{amsfonts}
\usepackage{amsmath}
\usepackage{amssymb}
\usepackage{amsthm}
\usepackage{hyperref}
\usepackage{enumerate}

\numberwithin{equation}{section}

\hypersetup{hidelinks}

\newtheorem{theorem}{Theorem}[section]
\newtheorem{proposition}[theorem]{Proposition}
\newtheorem{lemma}[theorem]{Lemma}
\newtheorem{remark}[theorem]{Remark}

\newtheorem{corollary}[theorem]{Corollary}

\newcommand{\pk}[2][\mathcal H]{\mathcal P_{#2} \left(#1\right)}
\newcommand{\fs}[1][\mathcal H]{\mathcal F_s \left(#1\right)}
\newcommand{\un}[1][k]{\mathcal U_{#1}}
\newcommand{\hk}[1][k]{\mathcal H_{#1}}
\newcommand{\mppq}[3][p]{\mathcal M_{#1} \left(#2, #3\right)}
\newcommand{\eppq}[3][p]{\mathcal E_{#1} \left(#2, #3\right)}
\newcommand{\NN}{\mathbb N}
\newcommand{\RR}{\mathbb R}
\newcommand{\CC}{\mathbb C}
\newcommand{\HH}{\mathcal H}
\newcommand{\KK}{\mathcal K}

\DeclareMathOperator{\tr}{tr}
\DeclareMathOperator{\rk}{rk}
\DeclareMathOperator{\ran}{ran}
\DeclareMathOperator{\diag}{diag}

\title{Schatten isometries of Grassmann spaces}
\author{Lucijan Plevnik}
\address{University of Ljubljana, Faculty of Mathematics and Physics, Jadranska ulica 19, Ljubljana, Slovenia}
\email{
	lucijan.plevnik@fmf.uni-lj.si}
\author{Jyoti Rani}
\address{Indian Institute of Science Education and Research (IISER) Mohali, Knowledge City, S.A.S Nagar, Punjab 140306}
\email{jyotirani@iisermohali.ac.in}
\date{}

\begin{document}

\begin{abstract}
We characterize all (not necessarily surjective) Schatten $p$-isometries of the Grassmann space of rank $k$ projections on a complex Hilbert space for every $p\in(0,\infty]$. We prove that every such isometry is induced by a linear or conjugate-linear isometry of the underlying Hilbert space $\HH$, with the additional possibility of composition with orthocomplementation when $\dim \HH = 2k > 2$.

Our result generalizes Wigner's theorem, which plays a central role in the mathematical foundations of quantum mechanics. It also extends previous results on Schatten $p$-isometries of Grassmann spaces, including characterizations of surjective isometries by Geh\'er and Šemrl \cite{gehse} for $p=\infty$ and by Qian et al. \cite{qian} for $0<p<\infty$ and that of not necessarily surjective isometries by Geh\'er \cite{geh-grass} for $p=2$.
\end{abstract}

\maketitle

\section{Introduction}

Let $\HH$ be a complex Hilbert space of dimension greater than $1$. We say that a bounded linear operator on $\HH$ is a projection if it is idempotent and self-adjoint. Let $k \in \NN$, $k < \dim \HH$. The set $\pk{k}$ of all projections of rank $k$ can be equipped with several metrics; for more information on this topic we refer to \cite{qiu}. In this paper, we will consider Schatten (quasi)norms on the real vector space $\fs$ of all finite rank bounded self-adjoint operators on $\HH$, i.e. for any $p \in (0,\infty)$ we have
\begin{eqnarray*}
	\|A\|_p = \left(\tr |A|^p\right)^{\frac{1}{p}}, & A \in \fs.
\end{eqnarray*}
Recall that $\|\cdot\|_p$ is a norm if and only if $p \ge 1$ and that for $A \in \fs$, $\lim\limits_{p \to \infty} \|A\|_p$ equals the operator norm of $A$, denoted by $\|A\|_\infty$. In this paper, we describe isometries on the set $\pk{k}$ with respect to $\|\cdot\|_p$ for any fixed $p \in (0,\infty]$. We remark that for $0<p<1$, preservation of $\|\cdot\|_p$-distances is equivalent to being an isometry for $d_p(P,Q)=\|P-Q\|_p^p$ which is a metric, as it satisfies the triangle inequality by Rotfel'd inequality \cite[Theorem IV. 2. 14]{bhat}.

Our main result reads as follows.
\begin{theorem}\label{main}
	Let $p \in (0,\infty]$, $k \in \NN$, $k < \dim \HH$, and let $\phi \colon \pk{k} \to \pk{k}$ be a map. Then we have
	\begin{eqnarray}\label{piso}
		\|\phi(P)-\phi(Q)\|_p = \|P-Q\|_p, & P, Q \in \pk{k},
	\end{eqnarray}
	if and only if there exists a linear or conjugate-linear isometry $U \colon \HH \to \HH$ such that we have either
	\begin{eqnarray*}
		\phi(P) = UPU^\ast, & P \in \pk{k},
	\end{eqnarray*}
	or $k > 1$, $\dim \HH = 2k$, and
	\begin{eqnarray*}
		\phi(P) = I-UPU^\ast, & P \in \pk{k}.
	\end{eqnarray*}
\end{theorem}

Our research is motivated by Wigner's theorem, a central result in the mathematical foundations of quantum mechanics, originally published in 1931 (see \cite{wig} for an English translation). It states that a bijective map $\phi \colon \pk{1} \to \pk{1}$ satisfying
	\begin{eqnarray}\label{trpq}
		\tr \phi(P) \phi(Q) = \tr PQ, & P, Q \in \pk{1},
	\end{eqnarray}
is induced by a unitary or an anti-unitary operator on $\HH$. This theorem has been generalized in several ways, see e.g. \cite{chev, molnar-book, pankov-book}. Among such generalizations are descriptions of isometries on $\pk{k}$ for $0 < p \le \infty$. Indeed, for $P, Q \in \pk{1}$ we have $\|P-Q\|_p = 2^{\frac{1}{p}} \sin \varphi$ and $\tr PQ = \cos^2 \varphi$, where $\varphi \in \left[0,\frac{\pi}{2}\right]$ is the angle between the one-dimensional ranges $\ran P$ and $\ran Q$ and we regard $\frac{1}{\infty}$ as $0$ throughout the paper. Hence,
\begin{eqnarray}\label{ppq}
	\|P-Q\|_p = 2^{\frac{1}{p}} \sqrt{1-\tr PQ}, & P, Q \in \pk{1}.
\end{eqnarray}
Therefore, a map $\phi \colon \pk{1} \to \pk{1}$ satisfies \eqref{trpq} if and only if it satisfies \eqref{piso} for $k=1$ and any fixed $p \in (0, \infty]$. Thus, a noteworthy generalization of Wigner's theorem is the theorem \cite[Theorem 1.1]{gehse} due to Geh\'er and Šemrl. It states that if a surjective map $\phi \colon \pk{k} \to \pk{k}$ satisfies \eqref{piso} for $p = \infty$, then $\phi$ is induced by a unitary or an anti-unitary operator $U \colon \HH \to \HH$, while it can be additionally composed with the map $P \mapsto I-P$, $P \in \pk{k}$, in the case when $\dim \HH = 2k > 2$. An analogous result was proven to hold in the case $0 < p < \infty$ in the context of von Neumann algebras \cite{qian}.

Next, note that in Wigner's theorem the assumption of injectivity is redundant since it easily follows from \eqref{trpq}. In fact, the surjectivity assumption can also be omitted, provided that the isometry $U$ is not required to be surjective; see, for example, \cite{geh-wig}. When $\HH$ is finite-dimensional, Gehér and Šemrl \cite{gehse} showed that isometries on $\pk{k}$ are automatically surjective for the operator norm. The same conclusion holds for every $p \in (0,\infty]$, as noted in Remark \ref{fin}, although our proof of the main theorem does not use this observation. In \cite{geh-grass}, Geh\'er characterized the not necessarily surjective isometries of $\pk{k}$ for $p=2$ on arbitrary complex Hilbert spaces. Our main result removes the assumption of surjectivity also in the infinite-dimensional case when $p \ne 2$. It thereby completes the classification of Schatten $p$-isometries of fixed finite-rank Grassmann spaces over complex Hilbert spaces for any $p \in (0,\infty]$.

In \cite{geh-grass}, Geh\'er expressed the hope that his proof for $p=2$ would provide insight into the nonsurjective operator-norm problem. In order to obtain a unified classification for all $p \in (0,\infty]$, we indeed build on the finite-dimensional reduction strategy of \cite{geh-grass}. When $\dim H>2k$, the key steps are to establish preservation of orthogonality and of the relation $R\leq P+Q$ for orthogonal $P, Q \in \pk{k}$. These properties allow us to reduce the problem to the case $\dim \HH=2k$, which is treated using preservation of adjacency in both directions for finite $p$ and the existing classification for the operator norm. The metric characterizations required for this reduction involve different arguments in the ranges $0<p<2$, $2\leq p<\infty$, and $p=\infty$.

The main additional difficulty occurs in the operator-norm case, where maximal distance does not characterize orthogonality when $k>1$. Following the approach from \cite{gehse}, we study equidistance sets associated with pairs of projections. However, the characterization used in \cite{gehse} does not directly yield preservation of orthogonality under nonsurjective isometries. For pairs at maximal distance, we instead characterize orthogonality by the existence of a scaled isometric copy of the unitary group $\un$, equipped with the operator-norm metric, in the associated equidistance set. This existence condition is preserved by isometries and therefore allows us to establish orthogonality preservation without assuming surjectivity.

\section{Preliminaries}

We denote the set of all projections of finite rank by $\pk{f}$. The following lemma is Wedin's version of the Halmos two projections theorem, see e.g. \cite[Corollary 2.2]{bot}.

\begin{lemma}\label{wed}
    Let $P, Q \in \pk{f}$ and denote $X = \ran P + \ran Q$. Then there exist $m, q, r, s \in \NN \cup \{0\}$, $\varphi_1, \ldots, \varphi_s \in \left(0, \tfrac{\pi}{2}\right)$, and an orthonormal basis of $X$, according to which we have matrix representations
	$$
		P|_X =
		I_m \oplus I_q \oplus 0_r \oplus
		\bigoplus_{j=1}^s \diag(1,0)
	$$
	and
	$$
		Q|_X =
		I_m \oplus 0_q \oplus I_r \oplus
		\bigoplus_{j=1}^s \begin{bmatrix}
		    \cos^2 \varphi_j & \cos \varphi_j \sin \varphi_j \\
			\cos \varphi_j \sin \varphi_j & \sin^2 \varphi_j
		\end{bmatrix}.
	$$
\end{lemma}

\begin{remark}
    We say that $\varphi_1, \ldots, \varphi_s$ from Lemma \ref{wed} are principal angles between $P$ and $Q$ from the interval $\left(0,\tfrac{\pi}{2}\right)$, while $0$ and $\tfrac{\pi}{2}$ are principal angles of respective multiplicities $m$ and $\min\{q,r\}$. In this way we get principal angles
    $$
    \varphi_1, \ldots, \varphi_k \in \left[0,\tfrac{\pi}{2}\right]
    $$
    when $\min\{\rk P, \rk Q\} = k$. An interested reader can find more about principal angles e.g. in \cite{bot}.
\end{remark}

The following result can be deduced directly from Lemma \ref{wed}.

\begin{corollary}\label{princdp}
	Let $P, Q \in \pk{k}$ with principal angles $0 \le \varphi_k \le \ldots \le \varphi_1 \le \tfrac{\pi}{2}$. Then
	$$
	\|P-Q\|_p = \left\{
    \begin{array}{ccc}
        \left(2\sum_{j=1}^k\sin^p\varphi_j\right)^{\frac{1}{p}} & : & 0 < p < \infty  \\
        \sin \varphi_1 & : & p = \infty
    \end{array}
    \right..
	$$
    If $0 < p < \infty$, then
    $$
    \|P-Q\|_p =
    \left(2\tr \left( P-PQP \right)^{\frac{p}{2}}\right)^{\frac{1}{p}}.
    $$
\end{corollary}

We say that $P, Q \in \pk{k}$ are {\em orthogonal} and write $P \perp Q$ if their ranges are orthogonal. This is equivalent to all principal angles being $\frac{\pi}{2}$, while an equivalent algebraic condition is $PQ=0$.

The following result is known in a more general context (\cite[Lemma 2.5]{qian}), but for our purposes it follows easily from Corollary \ref{princdp}.

\begin{corollary}\label{perp}
	Suppose that $\dim \HH \ge 2k$ and we have either $0 < p < \infty$ or $p = \infty$ and $k=1$. Then
	\begin{eqnarray*}
	P \perp Q \iff \|P-Q\|_p = (2k)^{\frac{1}{p}} & P, Q \in \pk{k}.
	\end{eqnarray*}
\end{corollary}

For projections $P$ and $Q$ we write $P \le Q$ if $\ran P \subseteq \ran Q$.

\begin{lemma}\label{pzpxpy}
    Let $P \in \pk{k}$ and let $Q$ be a projection. Then
    $$
    \tr PQP \le k
    $$
    and
    $$
    PQP = P \iff \tr PQP = k \iff P \le Q.
    $$
\end{lemma}

\begin{proof}
    We have $k - \tr(PQP) = \|P(I-Q)P\|_1 = \|(I-Q)P\|_2^2$.
\end{proof}

We say that $P, Q \in \pk{k}$ are {\em adjacent} and write
$$
P \sim Q
$$
if $\dim (\ran P \cap \ran Q) = k-1$.

\begin{lemma}\label{adjd}
    Suppose that $p \in (0,\infty) \setminus \{2\}$ and $\dim \HH = 2k$. Let $P, Q \in \pk{k}$. Then $P \sim Q$ if and only if
    \begin{equation}\label{zlt}
    0 < \|P-Q\|_p^p \le 2
    \end{equation}
    and
    \begin{equation}\label{pxybp}
    \frac{\|P-(I-Q)\|_p^p}{2} = k-1 + \left(1 - \left(\frac{\|P-Q\|_p^p}{2}\right)^{\frac{2}{p}}\right)^{\frac{p}{2}}.
    \end{equation}
\end{lemma}

\begin{proof}
    Let $0 \le \varphi_k \le \ldots \le \varphi_1 \le \tfrac{\pi}{2}$ be the principal angles between $P$ and $Q$. We obtain from Lemma \ref{wed} that
    \begin{eqnarray*}
    \frac{\|P-Q\|_p^p}{2} = \sum_{j=1}^k \sin^p \varphi_j & {\rm and} & \frac{\|P-(I-Q)\|_p^p}{2} = \sum_{j=1}^k \cos^p \varphi_j.
    \end{eqnarray*}
    We have
    $$
    P \sim Q \iff 0 < \varphi_1 \le \tfrac{\pi}{2} \textrm{ and } \varphi_2 = \ldots = \varphi_k = 0.
    $$
    Hence, $P \sim Q$ implies \eqref{zlt}.
    
    On the other hand, assume that \eqref{zlt} holds. Let $s = \frac{\|P-Q\|_p^p}{2}$. Then
    \begin{eqnarray*}
        s > 0 & {\rm and} & 0 < \varphi_1 \le \tfrac{\pi}{2}.
    \end{eqnarray*}
    Let $f_p \colon [0,1] \to \RR$ be given by $f_p(t) = \left(1-t^{\frac{2}{p}}\right)^{\frac{p}{2}}$. Then $f_p$ is a strictly convex function if $2 < p < \infty$, while it is strictly concave if $0<p<2$. Suppose first that $2 < p < \infty$. Then \eqref{zlt} implies that for all $j \in \{1, \ldots, k\}$ we have
    \begin{eqnarray*}
        0 \le \sin^p \varphi_j \le s \le 1
    \end{eqnarray*}
    and
    \begin{eqnarray*}
    \cos^p \varphi_j = f_p \left(\sin^p \varphi_j\right)
    \le \frac{\sin^p \varphi_j}{s} f_p \left(s\right) + \left(1-\frac{\sin^p \varphi_j}{s}\right) f_p(0) \\
    = 1 + \frac{\sin^p \varphi_j}{s} \left(f_p \left(s\right) - 1\right).
    \end{eqnarray*}
    By summing over all $j \in \{1, \ldots, k\}$ we obtain
    $$
    \frac{\|P-(I-Q)\|_p^p}{2} \le k-1 + f_p (s).
    $$
    If $0<p<2$, then all inequalities are reversed. Strict convexity/concavity of $f_p$ yields that in both cases, \eqref{pxybp} holds if and only if
    \begin{eqnarray*}
        \sin^p \varphi_j \in \left\{0,s\right\}, & j = 1, \ldots, k.
    \end{eqnarray*}
    Since $\sin \varphi_1 > 0$, the latter is equivalent to $\varphi_2 = \ldots = \varphi_k = 0$.
\end{proof}

\section{Containment for $0 < p < 2$}

For all $p \in (0,\infty)$ and $P, Q \in \pk{k}$ we denote
$$
\eppq{P}{Q} = \left\{ R \in \pk{k} \mid \|P-R\|_p^p + \|R-Q\|_p^p = 2k \right\}.
$$

\begin{lemma}\label{charzsxy}
    Suppose that $\dim \HH \ge 2k$ and $0<p<2$. Let $P, Q \in \pk{k}$ be orthogonal. Then
    $$
    \eppq{P}{Q} = \left\{ R \in \pk{k} \mid R \le P+Q, \, PR=RP \right\}.
    $$
\end{lemma}

\begin{proof}
    Let $R \in \pk{k}$. By Corollary \ref{princdp},
    $$
    \frac{\|P-R\|_p^p + \|R-Q\|_p^p}{2} = \tr \left( R-RPR \right)^{\frac{p}{2}} + \tr \left( R-RQR \right)^{\frac{p}{2}}.
    $$
    Since the spectra of $R-RPR$ and $R-RQR$ are contained in $[0,1]$, we deduce from $0 < p < 2$ that
    $$
    \tr \left( R-RPR \right)^{\frac{p}{2}} \ge \tr \left( R-RPR \right) = k - \tr \left(RPR\right)
    $$
    and
    $$
    \tr \left( R-RQR \right)^{\frac{p}{2}} \ge k - \tr \left(RQR\right),
    $$
    implying that
    $$
    \frac{\|P-R\|_p^p + \|R-Q\|_p^p}{2} \ge 2k - \tr \left(R\left(P+Q\right)R\right).
    $$
        It follows from $P \perp Q$ that $P+Q \in \pk{2k}$. Lemma \ref{pzpxpy} yields
        $$
        \tr \left(R\left(P+Q\right)R\right) \le k
        $$
        and
        $$
        \tr \left(R\left(P+Q\right)R\right) = k \iff R \le P+Q.
        $$
        Therefore, $R \in \eppq{P}{Q}$ if and only if $R \le P+Q$ and
    \begin{eqnarray}\label{dproj}
        R-RPR \textrm{ and } R-RQR \textrm{ are projections}.
    \end{eqnarray}
    Finally, if $R \le P+Q$, then Lemma \ref{pzpxpy} yields $R\left(P+Q\right)R = R$, meaning that \eqref{dproj} is equivalent to
    $$
    R-RPR \textrm{ is a projection}.
    $$
    It easily follows from Lemma \ref{wed} that the latter condition is fulfilled if and only if $PR=RP$.
\end{proof}

\begin{lemma}\label{compcon}
    Suppose that $k > 1$, $\dim \HH \ge 2k$, and $P, Q \in \pk{k}$. Then there exist $R, S \in \pk{k}$ such that $PR=RP$, $RS=SR$, and $SQ=QS$.
\end{lemma}

\begin{proof}
    Lemma \ref{wed} yields the existence of $\varphi_1, \ldots, \varphi_k \in \left[0, \tfrac{\pi}{2}\right]$, a $(2k)$-dimensional subspace $X \subseteq \HH$, and an orthonormal basis of $X$, according to which we have matrix representations
\begin{eqnarray*}
    P|_X = \bigoplus_{j=1}^k \diag(1,0)
	 & {\rm and} &
		Q|_X =
		\bigoplus_{j=1}^k \begin{bmatrix}
		    \cos^2 \varphi_j & \cos \varphi_j \sin \varphi_j \\
			\cos \varphi_j \sin \varphi_j & \sin^2 \varphi_j
		\end{bmatrix}.
\end{eqnarray*}
    We distinguish two cases. First, assume that $k$ is even, so $k = 2m$ for some $m \in \NN$. Let $R \in \pk{k}$ be such that $\ran R \subset X$ and
    $$
    R|_X = \bigoplus_{j=1}^m \left(I_2 \oplus 0_2\right).
    $$
    Then $R$ commutes with $P$ and $Q$, so we get the required result by taking $S = R$.

    Suppose, on the other hand, that $k \ge 3$ is odd, so $k = 2m+1$ for some $m \in \NN$. Let $R, S \in \pk{k}$ be such that $\ran R, \, \ran S \subset X$ and
    $$
    R|_X = \left(\bigoplus_{j=1}^m \left(I_2 \oplus 0_2\right)\right) \oplus \diag(1,0)
    $$
    and
    $$
    S|_X = \begin{bmatrix}
		    \cos^2 \varphi_1 & \cos \varphi_1 \sin \varphi_1 \\
			\cos \varphi_1 \sin \varphi_1 & \sin^2 \varphi_1
		\end{bmatrix} \oplus \bigoplus_{j=1}^m \left(I_2 \oplus 0_2\right).
    $$
    Then $R$ commutes with $P$ and $S$, while the latter commutes with $Q$.
\end{proof}

\begin{remark}
    If $k$ is even, then we may take $S = R$.
\end{remark}

\begin{lemma}\label{wow}
    Suppose that $k > 1$, $\dim \HH \ge 2k$, and $0 < p < 2$. Let $P, Q, R \in \pk{k}$ such that $P \perp Q$. Then $R \le P+Q$ if and only if there exist orthogonal $S, T \in \eppq{P}{Q}$ and orthogonal $U, V \in \eppq{S}{T}$ such that $R \in \eppq{U}{V}$.
\end{lemma}

\begin{proof}
    Assume first that there exist such $S, T, U, V$. Lemma \ref{charzsxy} implies that $S, T \le P+Q$, $U, V \le S+T$, and $R \le U+V$. Hence, $U+V \le S+T \le P+Q$. Since $U+V$, $S+T$, and $P+Q$ belong to $\pk{2k}$, they are equal. Consequently, $R \le P+Q$.

    Suppose, on the other hand, that $R \le P+Q$. Lemma \ref{compcon}, applied for the $(2k)$-dimensional Hilbert space $\ran P + \ran Q$, yields the existence of $S, U \in \pk{k}$ such that $S, U \le P+Q$, $PS=SP$, $SU=US$, and $UR=RU$. Let $T = P+Q-S$ and $V = P+Q-U$. Then $T, V \in \pk{k}$, $S \perp T$, $V \perp U$, and Lemma \ref{charzsxy} implies that $S, T \in \eppq{P}{Q}$, $U, V \in \eppq{S}{T}$, and $R \in \eppq{U}{V}$.
\end{proof}

\section{Containment for $2 \le p \le \infty$}

Following the approach from \cite{gehse}, for all $p \in (0,\infty]$ and $P, Q \in \pk{k}$ denote
$$
\mppq{P}{Q} = \left\{ R \in \pk{k} \mid \|P-R\|_p = \|R-Q\|_p = \frac{(2k)^{\frac{1}{p}}}{\sqrt{2}} \right\}
$$

\begin{lemma}\label{xyz}
    Suppose that $\dim \HH \ge 2k$ and we have either $k=1$ or $2 \le p \le \infty$. If $P, Q \in \pk{k}$ are orthogonal and $R \in \mppq{P}{Q}$, then $R \le P+Q$.
\end{lemma}

\begin{proof}
    In the case $p = \infty$, the lemma is a special case of \cite[Corollary 2.3]{gehse}. Therefore, we assume that $0<p<\infty$, although we could adapt our approach for the operator norm by taking the limit $p \to \infty$ in \eqref{jens}.
    
    Let $R \in \mppq{P}{Q}$. Corollary \ref{princdp} yields
    $$
    \frac{1}{2} = \left(\frac{\|R-P\|_p}{(2k)^{\frac{1}{p}}}\right)^2 = \left(\frac{\tr \left( R-RPR \right)^{\frac{p}{2}}}{k}\right)^{\frac{2}{p}}.
    $$
    We claim that
    \begin{eqnarray}\label{jens}
    \left(\frac{\tr \left( R-RPR \right)^{\frac{p}{2}}}{k}\right)^{\frac{2}{p}} \ge \frac{\tr \left( R-RPR \right)}{k}.
    \end{eqnarray}
    Indeed, $R-RPR$ may be seen as an operator on the $k$-dimensional space $\ran R$. Thus, we clearly have an equality in \eqref{jens} if $k=1$ or $p=2$, while in the case when $k > 1$ and $p > 2$, we obtain \eqref{jens} by applying Jensen's inequality.
    It now follows from \eqref{jens} that
    \begin{eqnarray*}
    \frac{1}{2} \ge \frac{\tr \left( R-RPR \right)}{k}
    = 1 - \frac{\tr \left(RPR\right)}{k},
    \end{eqnarray*}
    implying that
    $$
    \tr \left(RPR\right) \ge \frac{k}{2}.
    $$
    Since the assumptions on $P$ and $Q$ are symmetric, we analogously obtain
    $$
    \tr \left(RQR\right) \ge \frac{k}{2}.
    $$
    Consequently,
    $$
    \tr \left(R\left(P+Q\right)R\right) \ge k.
    $$
    It follows from $P \perp Q$ that $P+Q \in \pk{2k}$. Therefore, Lemma \ref{pzpxpy} yields $R \le P+Q$.
\end{proof}

\begin{lemma}\label{xyzuv}
    Suppose that $\dim \HH \ge 2k$ and we have either $k=1$ or $2 \le p \le \infty$. Let $P, Q, R \in \pk{k}$ such that $P \perp Q$. Then $R \le P+Q$ if and only if there exist orthogonal $S, T \in \mppq{P}{Q}$ such that $R \in \mppq{S}{T}$.
\end{lemma}

\begin{proof}
    Suppose first that such $S$ and $T$ exist. Lemma \ref{xyz} implies that $S \le P+Q$, $T \le P+Q$, and $R \le S+T$. Consequently, $S+T \le P+Q$. Since $S+T$ and $P+Q$ belong to $\pk{2k}$, they are the same. Hence $R \le P+Q$.

    Assume, on the other hand, that $R \le P+Q$. Denote $X = \ran P + \ran Q$. Then Lemma \ref{wed} yields the existence of $\varphi_1, \ldots, \varphi_k \in \left[0,\tfrac{\pi}{2}\right]$ and an orthonormal basis of $X$, according to which we have matrix representations
	$$
		P|_X =
		\bigoplus_{j=1}^k \diag(1,0)
	$$
	and
	$$
		R|_X =
		\bigoplus_{j=1}^k \begin{bmatrix}
		    \cos^2 \varphi_j & \cos \varphi_j \sin \varphi_j \\
			\cos \varphi_j \sin \varphi_j & \sin^2 \varphi_j
		\end{bmatrix}.
	$$
    It follows from $P \perp Q$ that
	$$
	Q|_X = \bigoplus_{j=1}^k \diag(0,1).
	$$
    Let $S, T \le P+Q$ be such that
    \begin{eqnarray*}
        S|_X =
		\bigoplus_{j=1}^k \begin{bmatrix}
		    \frac{1}{2} & -\frac{i}{2} \\
			\frac{i}{2} & \frac{1}{2}
		\end{bmatrix}
    & {\rm and } &
    T|_X =
		\bigoplus_{j=1}^k \begin{bmatrix}
		    \frac{1}{2} & \frac{i}{2} \\
			-\frac{i}{2} & \frac{1}{2}
		\end{bmatrix}.
    \end{eqnarray*}
    Then $S, T \in \mppq{P}{Q}$ are orthogonal and we have $R \in \mppq{S}{T}$.
\end{proof}

\section{Orthogonality for the operator norm}

\begin{lemma}\label{isop}
    Let $m \in \NN$ be such that there exist $P_0, P_1, \ldots, P_m \in \pk{k}$ such that
    \begin{eqnarray}
        \label{pij} \|P_i-P_j\|_\infty = 1, & 1 \le i < j \le m, \\
        \label{pzj} \|P_0-P_j\|_\infty \le \tfrac{1}{\sqrt{2}}, & 1 \le j \le m.
    \end{eqnarray}
    Then $m \le 2^k$.
\end{lemma}

\begin{proof}
Let $\KK$ be the sum of the ranges of $P_0, P_1, \ldots, P_m$. Let $\bigwedge\nolimits^k\mathcal{K}$ denote the $k$th antisymmetric tensor power. We refer to \cite[Section I.5]{bhat} for the definition of $\bigwedge\nolimits^k\mathcal{K}$ and its properties, including the equality
\begin{equation}\label{det}
\left\langle
x_1\wedge\cdots\wedge x_k,\,
y_1\wedge\cdots\wedge y_k
\right\rangle
=
\det\bigl[\langle x_a,y_b\rangle\bigr]_{a,b=1}^k.
\end{equation}
For each $j=0,\ldots,m$, choose an orthonormal basis $\{e_{j,1},\ldots,e_{j,k}\}$ of $\ran P_j$. We deduce from \eqref{det} that for
\[
\xi_j=e_{j,1}\wedge\cdots\wedge e_{j,k}
\in \bigwedge\nolimits^k\mathcal{K}
\]
we have $\|\xi_j\|^2=\det I_k=1$.

Let $0 \le i<j\le m$ and let $0 \le \varphi_k\le\cdots\le \varphi_1 \le \frac{\pi}{2}$ be the principal angles between $P_i$ and $P_j$. Then the singular values of $\bigl[\langle e_{i,a},e_{j,b}\rangle\bigr]_{a,b=1}^k$
are $\cos\varphi_1,\ldots,\cos\varphi_k$, so \eqref{det} yields
\[
|\langle\xi_i,\xi_j\rangle|^2
=\prod_{\ell=1}^k\cos^2\varphi_\ell.
\]
If $i \ge 1$, then \eqref{pij} yields $\varphi_1 = \frac{\pi}{2}$ by Corollary \ref{princdp}, implying that $\langle\xi_i,\xi_j\rangle = 0$. Thus $\xi_1,\ldots,\xi_m$ form an orthonormal family in $\bigwedge^k\mathcal{K}$. If $i=0$, then \eqref{pzj} implies that $\varphi_{\ell}\le\tfrac{\pi}{4}$, $\ell = 1, \ldots, k$. Consequently,
\[
|\langle\xi_0,\xi_j\rangle|^2
=\prod_{\ell=1}^k\cos^2\varphi_{\ell}
\ge 2^{-k}.
\]
Applying Bessel's inequality we obtain
\[
1=\|\xi_0\|^2
\ge \sum_{j=1}^m|\langle\xi_0,\xi_j\rangle|^2
\ge m\,2^{-k}.
\]
Therefore $m\le 2^k$.
\end{proof}

We denote the sets of all $k \times k$ unitary and Hermitian matrices by $\un$ and $\hk$, respectively. We equip $\un$ with the metric induced by the operator norm $\|\cdot\|_\infty$.

Recall that in a metric space $(M,d)$ we say that a subset $A \subseteq M$ {\em resolves} $M$ if for any distinct $x, y \in M$ there exists $a \in A$ such that $d(x,a) \ne d(y,a)$.

\begin{lemma}\label{resol}
    Every nonempty open subset of $\un$ resolves $\un$.
\end{lemma}

\begin{proof}
    Let $\mathcal O \subseteq \un$ be a nonempty open subset and let $U, V \in \un$ such that
    \begin{eqnarray}\label{equid}
        \|U-W\|_\infty = \|V-W\|_\infty, & W \in \mathcal O.
    \end{eqnarray}
    If $U \in \mathcal O$ or $V \in \mathcal O$, then \eqref{equid} clearly implies the required result $U = V$. Assume now that $\mathcal O$ contains neither $U$ nor $V$.

    We claim that
    \begin{eqnarray*}
    \mathcal R = \left\{U Z \mid Z \in \un, \, \textrm{the eigenvalue of $Z$ at maximal distance} \right. \\ \left. \textrm{from $1$ is unique, simple, and different from $-1$} \right\}
    \end{eqnarray*}
    is an open and dense subset of $\un$. Indeed, openness follows from continuity of the eigenvalues and the strict separation of the unique simple maximizing eigenvalue from the remaining spectrum. Density follows by diagonalizing a unitary matrix and perturbing its eigenvalue arguments so that they lie in $(-\pi,\pi)$ and have pairwise distinct absolute values. Thus, $\mathcal R \cap \mathcal O$ contains some $W$. Let $A \in \hk$ such that $\|A\|_\infty \le \pi$ and $W = U e^{iA}$. Denote $a = \|A\|_\infty$. Then $0 < a < \pi$ and $A$ has a unique eigenvalue of modulus $a$ and this eigenvalue is simple. For all $t \in [0,1]$ we have
    \begin{eqnarray}\label{len}
        \|U-Ue^{itA}\|_\infty = 2\sin \tfrac{t a}{2} & {\rm and} & \|Ue^{itA}-W\|_\infty = 2\sin \tfrac{(1-t) a}{2}.
    \end{eqnarray}
    Since $\mathcal O$ is open and contains $Ue^{iA}$, there exists $\tau \in (0,1)$ such that $Ue^{itA} \in \mathcal O$, $\tau \le t \le 1$. We deduce from \eqref{len} and \eqref{equid} that
    \begin{eqnarray}\label{vg}
        \|V-Ue^{itA}\|_\infty = 2\sin \tfrac{t a}{2}, & \tau \le t \le 1.
    \end{eqnarray}
    In particular,
    \begin{equation}\label{iv}
    \|I_k-V^\ast W\|_\infty = \|V-W\|_\infty = 2\sin \tfrac{a}{2}.
    \end{equation}

    Regard $\CC^k$ as a real Euclidean space and equip its unit sphere $S$ with the standard spherical distance
    $$
    d(x,y) = \arccos {\rm Re}\, y^\ast x = 2 \arcsin \frac{\|x-y\|}{2}.
    $$    
    By \eqref{iv}, there exist $x \in S$ and $s \in \{a,-a\}$ such that $V^\ast Wx = e^{is}x$. We deduce from $d(Wx,Vx) = a < \pi$ that
    \begin{eqnarray*}
        \gamma(t) = e^{its} Vx, & 0 \le t \le 1,
    \end{eqnarray*}
    is the unique spherical arc of length $a$ from $Vx$ to $Wx = e^{is} Vx$. But for every $t \in [\tau,1]$ we have
    \begin{eqnarray*}
    a = d(Wx,Vx) \le d(Vx,Ue^{itA}x) + d(Ue^{itA}x,Wx) \\
    = 2 \left(\arcsin \frac{\|Vx-Ue^{itA}x\|}{2} + \arcsin \frac{\|Ue^{itA}x-Wx\|}{2}\right) \\
    \le 2 \left(\arcsin \frac{\|V-Ue^{itA}\|_\infty}{2} + \arcsin \frac{\|Ue^{itA}-W\|_\infty}{2}\right) = a,
    \end{eqnarray*}
    where we applied \eqref{vg} and \eqref{len} in order to obtain the last equality. Hence,
    \begin{eqnarray*}
        d(Vx,Ue^{itA}x) = t \, d(Wx,Vx) & {\rm and} & d(Ue^{itA}x,Wx) = (1-t) d(Wx,Vx)
    \end{eqnarray*}
    for all $t \in [\tau,1]$, implying that
    \begin{eqnarray}\label{uei}
        Ue^{itA}x = \gamma(t), & \tau \le t \le 1.
    \end{eqnarray}
    Therefore,
    \begin{eqnarray*}
        e^{i(t-1)A}x = e^{i(t-1)s}x, & \tau \le t \le 1.
    \end{eqnarray*}
    By differentiating this equality at $t=1$ from the left we obtain $Ax=sx$. Hence, $e^{itA}x = e^{its} x$, so \eqref{uei} yields $Ux = Vx$. Since $s$ is the unique eigenvalue of $A$ of modulus $a$ and this eigenvalue is simple, its eigenspace is spanned by $x$.

    Complete $\{x\}$ to an orthonormal basis $\{x,x_2,\ldots,x_k\}$ of $\CC^k$. For each $j \in \{2, \ldots, k\}$ and $t \in [0,1]$, choose $T_j(t) \in \un$ such that
\[
T_j(t)x=\cos t\,x+\sin t\,x_j
\]
    and $\lim_{t \downarrow 0} \|T_j(t) - I_k\|_\infty = 0$. Then $UT_j(t)e^{iA}T_j(t)^\ast \in \mathcal R$ and since $\mathcal O$ is open, we have $UT_j(t)e^{iA}T_j(t)^\ast \in \mathcal O$ for all $t>0$ sufficiently close to $0$. The eigenspace of $T_j(t)AT_j(t)^\ast$ corresponding to its unique eigenvalue of modulus \(a\) is spanned by $T_j(t)x$. Hence the preceding argument yields $UT_j(t)x=VT_j(t)x$. Therefore, $Ux_j=Vx_j$, $j = 2, \ldots, k$, yielding $U=V$.
\end{proof}

\begin{lemma}\label{lorti}
    Suppose that $\dim \HH \ge 2k$ and let $P, Q \in \pk{k}$. Then $P \perp Q$ if and only if $\|P-Q\|_\infty = 1$ and there exists a map $\psi \colon \un \to \mppq[\infty]{P}{Q}$ such that
    \begin{eqnarray}\label{isu}
        \|\psi(U) - \psi(V)\|_\infty = \frac{\|U-V\|_\infty}{2}, & U, V \in \un.
    \end{eqnarray}
\end{lemma}

\begin{proof}
    By Corollary \ref{princdp}, $P \perp Q$ clearly implies $\|P-Q\|_\infty = 1$.

    Suppose now that $\|P-Q\|_\infty = 1$. If $P \perp Q$, then we have the orthogonal decomposition $\HH = \ran P \oplus \ran Q \oplus (\ran P \oplus \ran Q)^\bot$ and the map
    \begin{eqnarray*}
    U \mapsto \left(\tfrac{1}{2}\begin{bmatrix}
        I_k & U \\
        U^\ast & I_k
    \end{bmatrix} \right) \oplus 0, & U \in \un,
    \end{eqnarray*}
    clearly maps to $\mppq[\infty]{P}{Q}$ and satisfies \eqref{isu}.

    Assume, on the other hand, that there exists a map $\psi \colon \un \to \mppq[\infty]{P}{Q}$ satisfying \eqref{isu}. Seeking a contradiction, suppose that $P \not \perp Q$. It follows from \cite[Corollary 2.3]{gehse} that there exist $r \in \NN$, $r<k$, an orthogonal decomposition $\HH = X \oplus Y \oplus \KK$ such that $\dim X = \dim Y = r$, $P_0 \in \pk[\KK]{k-r}$, and maps $\alpha \colon \un \to \un[r]$ and $\beta \colon \un \to \pk[\KK]{k-r}$ such that
    \begin{eqnarray*}
        \psi(U) = \left(\tfrac{1}{2}\begin{bmatrix}
        I_r & \alpha(U) \\
        \alpha(U)^\ast & I_r
    \end{bmatrix}\right) \oplus \beta(U), & \|P_0 - \beta(U)\|_\infty \le \tfrac{1}{\sqrt{2}}, & U \in \un,
    \end{eqnarray*}
    and
    \begin{eqnarray*}
        \max \left\{ \|\alpha(U) - \alpha(V)\|_\infty, 2\|\beta(U) - \beta(V)\|_\infty \right\} = \|U-V\|_\infty, & U, V \in \un.
    \end{eqnarray*}
    In particular, for all $U, V \in \un$ we have
    \begin{eqnarray*}
        \|\alpha(U) - \alpha(V)\|_\infty \le \|U-V\|_\infty & {\rm and} & \|\beta(U) - \beta(V)\|_\infty \le \tfrac{\|U-V\|_\infty}{2}.
    \end{eqnarray*}
     Thus, $\alpha$ and $\beta$ are continuous.
    
    Let
    $$
    \mathcal S = \left\{ (U,V) \in \un \times \un \mid \|\beta(U) - \beta(V)\|_\infty < \tfrac{\|U-V\|_\infty}{2} \right\}.
    $$
    We distinguish two cases. First, suppose that $\mathcal S = \emptyset$. Then
    \begin{eqnarray*}
        \|\beta(U) - \beta(V)\|_\infty = \tfrac{\|U-V\|_\infty}{2}, & U, V \in \un.
    \end{eqnarray*}
    Hence, for each $\epsilon \in \{-1,1\}^k$ we have $\beta(\diag(\epsilon)) \in \pk[\KK]{k-r}$,
    $$
    \|P_0-\beta(\diag(\epsilon))\|_\infty \le \tfrac{1}{\sqrt{2}},
    $$
    and
    \begin{eqnarray*}
        \|\beta(\diag(\epsilon))-\beta(\diag(\delta))\|_\infty = 1, & \delta \in \{-1,1\}^k \setminus \{\epsilon\}.
    \end{eqnarray*}
    Lemma \ref{isop} yields $2^k \le 2^{k-r}$, a contradiction.
    
    Second, suppose that $\mathcal S \ne \emptyset$. Continuity of $\beta$ implies that $\mathcal S$ is open in $\un \times \un$. Hence, $\mathcal S$ contains $\mathcal U \times \mathcal V$ for some nonempty open subsets $\mathcal U, \mathcal V \subseteq \un$. We have
    \begin{eqnarray*}
        \|\alpha(U) - \alpha(V)\|_\infty = \|U-V\|_\infty, & U \in \mathcal U, \, V \in \mathcal V.
    \end{eqnarray*}
    We claim that $\alpha$ is injective on $\mathcal U$. Indeed, let $U_1, U_2 \in \mathcal U$ such that $\alpha(U_1) = \alpha(U_2)$. For all $V \in \mathcal V$ we have
    $$
    \|U_1-V\|_\infty = \|\alpha(U_1) - \alpha(V)\|_\infty = \|\alpha(U_2) - \alpha(V)\|_\infty = \|U_2-V\|_\infty.
    $$
    Lemma \ref{resol} yields $U_1=U_2$. Thus, $\alpha$ is an injective continuous map from an open subset of the $k^2$-dimensional real manifold $\un$ to the $r^2$-dimensional real manifold $\un[r]$, contradicting the invariance of domain theorem.
\end{proof}

\section{Proof of Theorem \ref{main}}

One implication is clear. So, let $p \in (0,\infty]$, $k \in \NN$, $k < \dim \HH$, and let $\phi \colon \pk{k} \to \pk{k}$ be a map satisfying \eqref{piso}.

\begin{lemma}\label{ort}
	If $\dim \HH \ge 2k$, then
	\begin{eqnarray*}
	P \perp Q \implies \phi(P) \perp \phi(Q), & P, Q \in \pk{k}.
	\end{eqnarray*}
\end{lemma}

\begin{proof}
    If either $0 < p < \infty$ or $k=1$, then we apply Corollary \ref{perp} and \eqref{piso}.

    Assume now that $p = \infty$ and let $P, Q \in \pk{k}$ be orthogonal. By Lemma \ref{lorti}, $\|P-Q\|_\infty = 1$ and there exists a map $\psi \colon \un \to \mppq[\infty]{P}{Q}$ such that
    \begin{eqnarray*}
        \|\psi(U) - \psi(V)\|_\infty = \frac{\|U-V\|_\infty}{2}, & U, V \in \un.
    \end{eqnarray*}
    It follows from \eqref{piso} that $\|\phi(P)-\phi(Q)\|_\infty = 1$ and
    $$
    \phi \circ \psi \colon \un \to \mppq[\infty]{\phi(P)}{\phi(Q)}
    $$
    satisfies
    \begin{eqnarray*}
        \|(\phi \circ \psi)(U) - (\phi \circ \psi)(V)\|_\infty = \frac{\|U-V\|_\infty}{2}, & U, V \in \un.
    \end{eqnarray*}
    Another application of Lemma \ref{lorti} yields $\phi(P) \perp \phi(Q)$.
\end{proof}

\begin{lemma}\label{adjp}
    If $k>1$, $\dim \HH = 2k$, and $0 < p < \infty$, then
	\begin{eqnarray*}
	P \sim Q \iff \phi(P) \sim \phi(Q), & P, Q \in \pk{k}.
	\end{eqnarray*}
\end{lemma}

\begin{proof}
    In the case $p = 2$, the lemma was proven in \cite[Subsection 2.1]{geh-grass}.
    
    Assume now that $p \in (0,\infty) \setminus \{2\}$. Since $\dim \HH = 2k$, Lemma \ref{ort} yields
    \begin{eqnarray}\label{fort}
        \phi(I-P) = I-\phi(P), & P \in \pk{k}.
    \end{eqnarray}
    We now apply Lemma \ref{adjd}, \eqref{piso}, and \eqref{fort}.
\end{proof}

\begin{proposition}\label{tk}
    If $k> 1$ and $\dim \HH = 2k$, then there exists a unitary or an anti-unitary operator $U \colon \HH \to \HH$ such that we have either
	\begin{eqnarray*}
        \phi(P) = UPU^\ast, & P \in \pk{k},
	\end{eqnarray*}
	or
	\begin{eqnarray*}
        \phi(P) = I-UPU^\ast, & P \in \pk{k}.
	\end{eqnarray*}
\end{proposition}

\begin{proof}
    In the case $p = \infty$, the proposition was proven in \cite{gehse}.
    
    If $0 < p < \infty$, then Lemma \ref{ort}, Lemma \ref{adjp}, and $k > 1$ imply that we can apply a nonsurjective version of Chow's theorem for orthogonality preserving maps (see \cite[Subsection 5.3]{pankov}) which yields the required result.
\end{proof}

\begin{remark}\label{fin}
    For \(0<p<\infty\), we can alternatively prove Proposition \ref{tk} by showing that $\phi$ is surjective and applying \cite[Theorem 3.4]{qian}.

    Automatic surjectivity of $\phi$ can be established under the assumption $\dim \HH < \infty$. Indeed, for $\alpha=\min\{1,p\}$,
\[
d(P,Q)=\|P-Q\|_p^\alpha,
\qquad P,Q\in\pk{k},
\]
is a metric which induces the operator norm topology on $\pk{k}$, because
\[
\|P-Q\|_\infty^\alpha
\le d(P,Q)
\le (\dim\HH)^{\alpha/p}\|P-Q\|_\infty^\alpha.
\]
Hence, $\pk{k}$ is compact with respect to $d$, so surjectivity of $\phi$ follows from \cite[Theorem~1.6.14]{bur}.

Note also that automatic surjectivity in the finite-dimensional case was established for $p = \infty$ in \cite{gehse} by a different argument.
\end{remark}

\begin{lemma}\label{kopt}
    Suppose that $\dim \HH > 2k$. Let $P, Q, R \in \pk{k}$ such that $P \perp Q$ and $R \le P+Q$. Then $\phi(P) \perp \phi(Q)$ and $\phi(R) \le \phi(P)+\phi(Q)$.
\end{lemma}

\begin{proof}
    We apply \eqref{piso}, Lemma \ref{ort}, and either Lemma \ref{xyzuv} (when $k=1$ or $2 \le p \le \infty$) or Lemma \ref{wow} (when $k>1$ and $0 < p < 2$).
\end{proof}

\begin{proposition}\label{pgt}
    If $k>1$ and $\dim \HH > 2k$, then there exists a linear or conjugate-linear isometry $U \colon \HH \to \HH$ such that
	\begin{eqnarray*}
		\phi(P) = UPU^\ast, & P \in \pk{k}.
	\end{eqnarray*}
\end{proposition}

\begin{proof}
    Let $R \in \pk{2k}$ and $P_0 \in \pk{k}$ such that $P_0 \le R$ and let $Q_0 = R-P_0$. Then $P_0 \perp Q_0$ and $R = P_0+Q_0$. Denote $S = \phi(P_0)+\phi(Q_0)$. It follows from Lemma \ref{kopt} that $S \in \pk{2k}$ and
    \begin{eqnarray*}
        P \le R \implies \phi(P) \le S, & P \in \pk{k}.
    \end{eqnarray*}
    Thus, $\phi$ induces a map $\pk[\ran R]{k} \to \pk[\ran S]{k}$ which is an isometry with respect to $\|\cdot\|_p$. Proposition \ref{tk} yields the existence of a linear or conjugate-linear isometry $U_R \colon \ran R \to \HH$ such that $\ran U_R = \ran S$ and we have either
	\begin{eqnarray*}
		\phi(P) = U_RP|_{\ran R}U_R^\ast, & P \in \pk{k}, \, P \le R,
	\end{eqnarray*}
	or
	\begin{eqnarray*}
		\phi(P) = U_R(I-P)|_{\ran R}U_R^\ast, & P \in \pk{k}, \, P \le R.
	\end{eqnarray*}
    Since every pair $P, Q \in \pk{k}$ is dominated by some $R \in \pk{2k}$, we have
	\begin{eqnarray*}
		\tr \phi(P) \phi(Q) = \tr PQ, & P, Q \in \pk{k}.
	\end{eqnarray*}
    Hence, $\phi$ satisfies the assumptions of the main theorem in \cite{geh-grass} and we can obtain the required result from that theorem, although it suffices to apply Lemma 1 and the arguments from Section 2.2 of that paper in order to complete the proof.
\end{proof}

Proposition \ref{tk} and Proposition \ref{pgt} prove Theorem \ref{main} in the case $k>1$, $\dim \HH \ge 2k$.

In the case $k = 1$, our main tools Lemma \ref{ort} and Lemma \ref{kopt} work as well and we could easily adapt the rest of the proof in this case, but for the sake of simplicity we will rather cite the non-surjective version of Wigner's theorem \cite{geh-wig}.

In the remaining case $k < \dim \HH < 2k$ we can use a standard trick to apply the already established results for the isometry $Q \mapsto I - \phi(I-Q)$ on $\pk{\dim \HH-k}$.

%

\bibliographystyle{amsplain}
\bibliography{References}

@book{bhat,
 author = {Bhatia, Rajendra},
 title = {Matrix analysis},
 fseries = {Graduate Texts in Mathematics},
 series = {Grad. Texts Math.},
 issn = {0072-5285},
 volume = {169},
 isbn = {0-387-94846-5},
 year = {1996},
 publisher = {New York, NY: Springer},
 zbMATH = {967931},
 Zbl = {0863.15001}
}

@article{bot,
 author = {B{\"o}ttcher, A. and Spitkovsky, I. M.},
 title = {A gentle guide to the basics of two projections theory},
 fjournal = {Linear Algebra and its Applications},
 journal = {Linear Algebra Appl.},
 issn = {0024-3795},
 volume = {432},
 number = {6},
 pages = {1412--1459},
 year = {2010},
 doi = {10.1016/j.laa.2009.11.002},
 zbMATH = {5667879},
 Zbl = {1189.47073}
}

@book{bur,
 author = {Burago, D. and Burago, Yu. and Ivanov, S.},
 title = {A course in metric geometry},
 fseries = {Graduate Studies in Mathematics},
 series = {Grad. Stud. Math.},
 issn = {1065-7339},
 volume = {33},
 isbn = {0-8218-2129-6},
 year = {2001},
 publisher = {Providence, RI: American Mathematical Society (AMS)},
 zbMATH = {1626771},
 Zbl = {0981.51016}
}

@incollection{chev,
 author = {Chevalier, Georges},
 title = {Wigner's theorem and its generalizations},
 booktitle = {Handbook of quantum logic and quantum structures},
 isbn = {978-0-444-52870-4},
 pages = {429--475},
 year = {2007},
 publisher = {Amsterdam: Elsevier},
 zbMATH = {5214897},
 Zbl = {1126.81031}
}

@article{geh-wig,
 author = {Geh{\'e}r, Gy{\"o}rgy P{\'a}l},
 title = {An elementary proof for the non-bijective version of {Wigner}'s theorem},
 fjournal = {Physics Letters. A},
 journal = {Phys. Lett., A},
 issn = {0375-9601},
 volume = {378},
 number = {30-31},
 pages = {2054--2057},
 year = {2014},
 doi = {10.1016/j.physleta.2014.05.039},
 zbMATH = {6547005},
 Zbl = {1331.46066}
}

@article{geh-grass,
 author = {Geh{\'e}r, Gy{\"o}rgy P{\'a}l},
 title = {Wigner's theorem on {Grassmann} spaces},
 fjournal = {Journal of Functional Analysis},
 journal = {J. Funct. Anal.},
 issn = {0022-1236},
 volume = {273},
 number = {9},
 pages = {2994--3001},
 year = {2017},
 doi = {10.1016/j.jfa.2017.06.011},
 url = {real.mtak.hu/166882/1/1706.02329.pdf},
 zbMATH = {6774165},
 Zbl = {1372.81088}
}

@article{gehse,
 author = {Geh{\'e}r, Gy{\"o}rgy P{\'a}l and {\v{S}}emrl, Peter},
 title = {Isometries of {Grassmann} spaces},
 fjournal = {Journal of Functional Analysis},
 journal = {J. Funct. Anal.},
 issn = {0022-1236},
 volume = {270},
 number = {4},
 pages = {1585--1601},
 year = {2016},
 doi = {10.1016/j.jfa.2015.11.018},
 zbMATH = {6535753},
 Zbl = {1330.47045}
}

@book{molnar-book,
 author = {Moln{\'a}r, Lajos},
 title = {Selected preserver problems on algebraic structures of linear operators and on function spaces},
 fseries = {Lecture Notes in Mathematics},
 series = {Lect. Notes Math.},
 issn = {0075-8434},
 volume = {1895},
 isbn = {3-540-39944-5},
 year = {2007},
 publisher = {Berlin: Springer},
 doi = {10.1007/3-540-39944-5},
 zbMATH = {5076921},
 Zbl = {1119.47001}
}

@article{pankov,
 author = {Pankov, Mark},
 title = {Geometric version of {Wigner}'s theorem for {Hilbert} {Grassmannians}},
 fjournal = {Journal of Mathematical Analysis and Applications},
 journal = {J. Math. Anal. Appl.},
 issn = {0022-247X},
 volume = {459},
 number = {1},
 pages = {135--144},
 year = {2018},
 doi = {10.1016/j.jmaa.2017.10.065},
 zbMATH = {6817535},
 Zbl = {1396.46022}
}

@book{pankov-book,
 author = {Pankov, Mark},
 title = {Wigner-type theorems for {Hilbert} {Grassmannians}},
 fseries = {London Mathematical Society Lecture Note Series},
 series = {Lond. Math. Soc. Lect. Note Ser.},
 issn = {0076-0552},
 volume = {460},
 isbn = {978-1-108-79091-8; 978-1-108-80032-7},
 year = {2020},
 publisher = {Cambridge: Cambridge University Press},
 doi = {10.1017/9781108800327},
 zbMATH = {7144008},
 Zbl = {1466.81002}
}

@article{qian,
 author = {Qian, Wenhua and Shen, Junhao and Shi, Weijuan and Wu, Wenming and Yuan, Wei},
 title = {Surjective {{\(L^p\)}}-isometries on {Grassmann} spaces},
 fjournal = {Science China. Mathematics},
 journal = {Sci. China, Math.},
 issn = {1674-7283},
 volume = {66},
 number = {9},
 pages = {2105--2118},
 year = {2023},
 doi = {10.1007/s11425-021-2055-5},
 zbMATH = {7732872}
}

@article{qiu,
 author = {Qiu, Li and Zhang, Yanxia and Li, Chi-Kwong},
 title = {Unitarily invariant metrics on the {Grassmann} space},
 fjournal = {SIAM Journal on Matrix Analysis and Applications},
 journal = {SIAM J. Matrix Anal. Appl.},
 issn = {0895-4798},
 volume = {27},
 number = {2},
 pages = {507--531},
 year = {2005},
 doi = {10.1137/040607605},
 url = {repository.ust.hk/ir/bitstream/1783.1-1396/1/QZL.pdf},
 zbMATH = {5029248},
 Zbl = {1099.15024}
}

@misc{wig,
 author = {Wigner, Eugene P.},
 title = {Group theory and its application to the quantum mechanics of atomic spectra. {With} additions and corrections by {Eugene} {P}. {Wigner}},
 year = {1959},
 howpublished = {Pure and {Applied} {Physics}, 5. {New} {York}: {Academic} {Press} {Inc}., {Publ}.; {London}: {Academic} {Books} {Ltd}. {XI}, 372},
 zbMATH = {3140172},
 Zbl = {0085.37905}
}

\end{document}